\documentclass[oneside, a4paper,12pt]{article}
\usepackage{array,latexsym}
\usepackage[dvips]{graphicx}
\usepackage{amsfonts}
\usepackage{amsmath}
\usepackage{amsthm}
\usepackage{a4wide}
\usepackage{enumerate}
\newtheorem{theorem}{Theorem}
\newtheorem{lemma}[theorem]{Lemma}

\newtheorem{conjecture}[theorem]{Conjecture}
\theoremstyle{definition}
\newtheorem{definition}{Definition}
\theoremstyle{remark}

\newtheorem{example}{Example}
\newcommand{\pat}{2--13~}
\newcommand{\eqd}{${\sim}$~}
\newcommand{\eqc}{$\overset{c}{\sim}$~}
\newcommand{\qh}[1]{[#1]_{q}}
\newcommand{\qxh}[1]{[#1]^{x}_{q}}
\newcommand{\qyh}[1]{[#1]^{y}_{q}}
\newcommand{\qxyh}[1]{[#1]^{x,y}_{q}}

\newcommand{\MAP}{\mathbb{M}}
\newcommand{\BIJ}{\mathbb{B}}
\newcommand{\BIJB}{\mathbb{C}}

\begin{document}
\title{Permutations, cycles, and the pattern \pat}
\author{Robert Parviainen\\
  {\small ARC Centre of Excellence for Mathematics 
    and Statistics of Complex Systems}\\
 {\small 139 Barry Street, The University of Melbourne, Victoria, 3010}\\
 {\small E-mail: {\tt robertp@ms.unimelb.edu.au}}} 
\date{\small July 31, 2006\\
	%\dateline{Jan 1, 2006}{Xxx, 200x}\\
   }
\maketitle
\begin{abstract}
We count the number of occurrences of restricted patterns of length 3 in permutations with respect to length \emph{and} the number of cycles. The main tool is a bijection between permutations in standard cycle form and weighted Motzkin paths.
\end{abstract}
\section{Introduction}

Let $\mathcal{S}_{n}$ denote the set of permutations of $[n]=\{1,2,\dots, n\}$. A \emph{pattern} in a permutation $\pi\in \mathcal{S}_{n}$ is a permutation $\sigma\in\mathcal{S}_{k}$ and an occurrence of $\sigma$ as a subword of $\pi$:
There should exist $i_{1}<\dots < i_{k}$ such  that $\sigma=R(\pi(i_{1})\cdots\pi({i_{k}}))$, where $R$ is the reduction operator that maps the smallest element of the subword to 1, the second smallest to 2, and so on.

For example, an occurrence of the pattern 3--2--1 in $\pi\in \mathcal{S}_{n}$ means that there exists $1\leq i<j<k\leq n$ such that $\pi({i})>\pi({j})>\pi({k})$.

We further consider \emph{restricted} patterns, introduced by Babson and Steingr{\'\i}msson, \cite{BS2000}. The restriction is that two specified adjacent elements in the pattern \emph{must be adjacent} in the permutation as well. The position of the restriction in the pattern is indicated by \emph{an absence}  of a dash (--).  Thus,  an occurrence of the pattern 3--21 in $\pi\in \mathcal{S}_{n}$ means that there exists $1\leq i<j<n$ such that $\pi({i})>\pi({j})>\pi({j+1})$.

Here we are mainly interested in patterns of the type 2--13. We remark that it is shown by Claesson, \cite{Clas2001}, that the occurrences of \pat are equidistributed with the occurrences of the pattern 2--31, as well as  with 13--2 and with 31--2. The number of permutations with $k$ occurrences of \pat where given by Claesson and Mansour, \cite{CM2002}, for $k\leq 3$  and for $k\leq 8$ by Parviainen, \cite{Parv2006}.

The starting point of \cite{Parv2006} and this paper is a generating function related to the solution of a certain much studied Markov chain, the asymmetric exclusion process, \cite{BCEPR2006}. This function, of 4 variables, is the continued fraction
\begin{equation}
	F(q,x,y,t)=\\
	\cfrac{1}{1-t (\qxh{1}+\qyh{1})
   	-\cfrac{t^{2}\qh{1}\qxyh{2}}{1-t (\qxh{2}+\qyh{2})
	-\cfrac{t^{2}\qh{2}\qxyh{3}}{1-t (\qxh{3}+\qyh{3}) \dotsb}}},
\end{equation}
where
\begin{align*}
\qh{h}&=1+q+\cdots + q^{h-2} + q^{h-1},\\
\qxh{h}&=1+q+\cdots + q^{h-2}+ {x} q^{h-1},\\
\qyh{h}&=1+q+\cdots + q^{h-2}+ {y} q^{h-1},\\
\qxyh{h}&=1+q+\cdots + q^{h-3}+ ({x+y-xy}) q^{h-2}+{xy}q^{h-1}.
\end{align*}
It was shown in \cite{CM2002} and \cite{Parv2006} that $F(q,1,1,t)$ counts the number of permutations with $k$ occurrences of the pattern 2--13. The main goal of this paper is to study $F(q,x,1,t)$ and give a combinatorial interpretation of the coefficients. It turns out that the variable $x$ is connected to the cycle structure of permutations.

\section{Introducing cycles}
First consider $F(1,x,1,t)$, and expand in $t$:
\[
F(1,x,1,t)=1+(1+x)t+(2+3x+x^{2})t^{2}+(6+11x+6x^{2}+x^{3})t^{3}+O(t^{4}).
\]
These coefficients certainly looks like the unsigned Stirling numbers of the first kind. Thus $F(1,x,1,t)$ should count the number of permutations with respect to length and number of cycles.  This will indeed follow from the main theorem.

As $F(q,1,1,t)$ counts the number of occurrences of the pattern 2--13 and $F(1,x,1,t)$ the number of cycles, $F(q,x,1,t)$ should give (some kind of) bivariate statistic of occurrences of \pat and cycle distribution. 

\subsection{Cyclic occurrence of patterns}
The \emph{standard cycle form} of a permutation $\pi \in \mathcal{S}_{n}$ is the permutation written in cycle form, with cycles starting with the smallest element, and cycles ordered in decreasing order with respect to their minimal elements.  Let $C(\pi)$ denote the standard cycle form of a permutation $\pi$.
\begin{example}
If $\pi=47613852$, then $C(\pi)=(275368)(14)$.
\end{example}

\begin{definition}
Let $\pi$ be a permutation of $[n]$, with standard cycle form
\[
C(\pi)=(c_{1}^{1} c_{2}^{1}\cdots c_{i_{1}}^{1})(c_{1}^{2} c_{2}^{2}\cdots c_{i_{2}}^{2})\cdots (c_{1}^{k} c_{2}^{k}\cdots c_{i_{k}}^{k}),
\]
and let $\sigma=AB$ be a permutation of $[m]$, $m<n$.  The pattern  {\it A--B} occurs \emph{cyclically} in $\pi$ if occurs in one of the following senses
\begin{description}
	\item [Between cycles:]
		If {\it A--B} occurs in the permutation
		\[\hat \pi = c_{1}^{1} c_{2}^{1}\cdots c_{i_{1}}^{1}c_{1}^{2} c_{2}^{2}\cdots c_{i_{2}}^{2}\cdots c_{1}^{k} c_{2}^{k}\cdots c_{i_{k}}^{k}\] 
		\emph{and} there exists $a<b$ such that $A$ occurs in $c_{1}^{a} \cdots c_{i_{a}}^{a}$ and $B$ occurs in $c_{1}^{b} \cdots c_{i_{b}}^{b}$, we say that \pat occurs \emph{between cycles in $\pi$.}
	\item [Within cycles:] Let $\tilde \pi=c_{1}^{a}\cdots c_{i_{a}}^{a}$. If {\it A--B} occurs \ in $\tilde \pi$ we say that {\it A--B} occurs \emph{within cycle $a$ in $\pi$. }
\end{description}
\end{definition}
\begin{example}
If $C(\pi)=(275368)(14)$  there are 2 occurrences of 2--13 \emph{between} cycles, 2--14 and 3--14, and 2 occurrences of 2--13 \emph{within} cycles, 7--68 and 5--36.
\end{example}

Let $\Phi_{i,j}(n)$ denote the number of permutations of length $n$, with $i$ cyclic occurrences of \pat and $j$ cycles. 
\begin{theorem} \label{theorem:main} 
The function $F(q, x, 1, z)$ is the (ordinary) generating function for $\Phi_{i,j}(n)$:
	\[\Phi_{i,j}(n)=[q^{i}x^{j}z^{n}]F(q,x,1,z).\]
\end{theorem}  

%%%
%%%
%%%

\section{Proof of Theorem \ref{theorem:main}}
We will use the fact \cite[Theorem 1]{Flaj1980} that $F(q,x,1,t)$ is the generating function for weighted bi-coloured Motzkin paths.

\begin{definition}
  A \emph{Motzkin path} of length $n$ is a sequence of  vertices $p = (v_0,v_1,\dots,v_n)$, with $v_i \in \mathbb{N}^2$  (where $\mathbb{N} = \{0,1,\dots\}$), with steps $v_{i+1}-v_i \in \{ (1,1), (1,-1), (1,0)\}$ and $v_0 = (0,0)$ and $v_n=(n,0)$. 
  
   A  \emph{bicoloured Motzkin path} is a Motzkin path in which each east, $(1,0)$, 
  step is labelled by one of two colours. 
\end{definition}

From now on \emph{all Motzkin paths considered will be bi-coloured}. 

Let $N$ ($S$) denote a north, $(1,1)$, step (resp., south, $(1,-1)$, step), and $E$ and $F$ the two different coloured  east steps. Further, let $N_{h}, S_{h}, E_{h}, F_{h}$ denote the weight of a $N$, $S$, $E$, $F$ step, respectively, that starts at height $h$. The weight of a Motzkin path is the product of the steps weights. 

If the weights are given by
\begin{equation}\label{eq:weights 1}
	N_{h}=\qxh{h+2},  S_{h}=\qh{h}, E_{h}=\qxh{h+1} \mbox{ and } F_{h}=\qh{h+1},
\end{equation}
it follows immediately from \cite[Theorem 1]{Flaj1980} that $[q^{i}x^{j}t^{n}]F(q,x,1,t)$ is the number of Motzkin paths of length $n$ with weight $q^{i}x^{j}$. Let $\mathcal{M}_{n}$ denote the set of weighted Motzkin paths of length $n$ with step weights given by \eqref{eq:weights 1}.

To establish Theorem \ref{theorem:main} we will use a bijection between  sets of permutations  and weighted Motzkin paths of length $n$. 

%%%
%%%
%%%

\subsection{The arc representation}
\begin{figure}
	\begin{center}
		\includegraphics[scale=1.0]{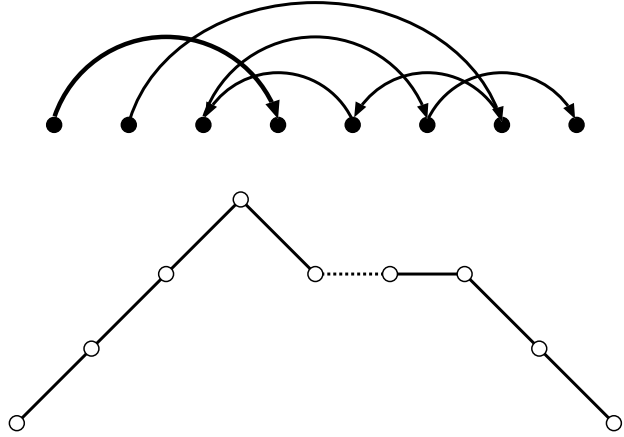}
		\caption{On top, the arc representation of $C(\pi)=(275368)(14)$. The node weights are, in order, $x, xq, q, 1, q, 1, q, 1$.  The shape pairs are, in order,  $(\emptyset,~\rightarrow), (\emptyset,~\rightarrow),  (\emptyset,~\rightleftharpoons), (\rightarrow,~\emptyset), (\leftarrow,~\leftarrow),  (\rightarrow,~\rightarrow), (\rightleftharpoons,~\emptyset), (\rightarrow,~\emptyset)$. At bottom, the image Motzkin path.}
		\label{fig:arcex}
	\end{center}
\end{figure}

We use a graphical representation of permutations to aid in the description of the mapping. For permutation $\pi\in\mathcal{S}_n$ with standard cycle form 
\[
C(\pi)=(c_{1}^{1} c_{2}^{1}\cdots c_{i_{1}}^{1})(c_{1}^{2} c_{2}^{2}\cdots c_{i_{2}}^{2})\cdots (c_{1}^{k} c_{2}^{k}\cdots c_{i_{k}}^{k}),
\]
make $n$ nodes in a line, representing the elements $1$ to $n$. For $s=1,\ldots k$ and $t=1,\ldots, i_{s}-1$ draw an arc from node $c_{t}^{s}$ to node $c_{t+1}^{s}$.   If cycle $s$ is of size 1 draw a loop from $c_{1}^{s}$ to itself. See Figure \ref{fig:arcex} for an example.

Associate each node with a left and a right shape. The left (right) shape is connections to nodes on the left (right) side with the node. The possible shapes on both sides are $\{\emptyset, \rightarrow, \leftarrow, \rightleftharpoons\}$. See Figure \ref{fig:arcex} for an example.

%%%
%%%
%%%

\subsubsection{Weights in the arc representation}
We now give each element, or node in the arc representation, a weight $x^{a}q^{b+w}$, in such a way that the product of a permutation's elements weights is $x^{k}q^{m}$, where $k$ is the number of cycles in $\pi$ and $m$ is the number of cyclic occurrences of 2--13. 

Imagine the arcs being drawn in sequence, in the order $c_{1}^{1}\rightarrow c_{2}^{1},  c_{2}^{1}\rightarrow c_{3}^{1}, \ldots , c_{i_{k}-1}^{k}\rightarrow c_{i_{k}}^{k}$. (If $i_{s}=1$ for a cycle $s$, we draw the loop $c_{1}^{s}\rightarrow c_{1}^{s}$.)  

Give node $k$ weight $x^{a}q^{b+w}$, where 
\begin{itemize}
	\item $a$ is 1 if the left right shape pair of the node is $(\emptyset, \rightarrow)$ and $a$ is 0 otherwise (element $k$ is the first in the cycle),
	\item $b$ is the number of times an arc \emph{belonging to a different cycle} that is drawn \emph{after} the node is visited passes over the node \emph{from left to right} (element $k$ plays the role of ``2'' in $b$ occurrences of 2--13 \emph{between} cycles),
	\item $w$ is the number of times an arc \emph{belonging to the same cycle} that is drawn \emph{after} the node is visited passes over the node \emph{from left to right} (element $k$ plays the role of ``2'' in $w$ occurrences of 2--13 \emph{within} cycles).
\end{itemize}
See Figure \ref{fig:arcex} for an example.

%%%
%%%
%%%

\subsection{One surjection and two bijections}
First we define a mapping $\MAP$ from $\mathcal{S}_{n}$ to Motzkin paths of length $n$, and prove that it is a surjection.

\begin{definition}
If $\pi\in\mathcal{S}_{n}$ have left shapes $\{l_{1},\ldots, l_{n}\}$ and right shapes $\{r_{1},\ldots, r_{n}\}$, let step $k$ in $\MAP(\pi)$ be $s_{k}$, where $s_{k}$ is given by the following table (where ``$-$'' denotes pairs of shapes that do not appear). Further, give step $k$ the same weight as node $k$.
\begin{center}
	\begin{tabular}{c||c|c|c|c}
		$l_{k}\backslash r_{k}$ & $\emptyset$ & $\rightarrow$ & $\leftarrow$ & $\rightleftharpoons$ \\\hline\hline 
		$\emptyset$ & $E$ & $N$ & $F$ & $N$ \\\hline 
		$\rightarrow$ & $S$ & $E$ & $-$ & $-$ \\\hline 
		$\leftarrow$ & $-$ & $-$ & $F$ & $-$ \\\hline 
		$\rightleftharpoons$ & $S$ & $-$ & $-$ & $-$
	\end{tabular}
\end{center}
\end{definition}
\begin{lemma}
The the mapping  $\MAP$ is a surjection from the set of permutations to the set of Motzkin paths with no $F$ steps at level 0.
\end{lemma}
\begin{proof}
To show that the image is a Motzkin path, the conditions for a Motzkin path must be verified. Namely, that $N_{k}\geq S_{k}$, $k<n$, $N_{n}=S_{n}$ where $N_{m}$ and $S_{m}$ are the number of $N$ and $S$ steps, respectively, up to and including step $m$.

As the shape pairs $(\leftarrow, \leftarrow)$,  $(\rightarrow, \rightarrow)$, $(\emptyset, \emptyset)$ and $(\emptyset, \leftarrow)$ map to $E$ and $F$ steps, it may be assumed that these shapes do not occur.

Now, in a valid arc diagram, the number of $(\emptyset, \rightarrow)$ and $(\emptyset, \rightleftharpoons)$ shape pairs up to and including node $k$ must be greater than or equal to the number of $(\rightarrow, \emptyset)$ and $(\rightleftharpoons, \emptyset)$ shape pairs. Further, these counts must agree for $k=n$. This is exactly what is needed. 

To show that $\MAP$ is a surjection, consider any Motzkin path $p$ with no $F$ steps at level 0. We will build an arc diagram $a$ that maps to $p$. 

For each $F$ step in $p$ we can associate a unique pair of $N$ and $S$ steps. (The rightmost (leftmost) $N$ ($S$) step to the left (right) of the $F$ step, ending (starting) at the $F$ steps level.) Let $n, f, s$ denote the positions of the $N, F, S$ steps, respectively. In the arc diagram $a$, draw an arc from node $n$ to node $s$ and one from node $s$ to node $f$. 

For the remaining $N$ and $S$ steps, fix one pairing of these, and draw arcs from the nodes corresponding to the $N$ steps, to the associated nodes corresponding to the $S$ steps.
 
For every $E$ step in $p$, draw an loop at the corresponding node.

Clearly, $a$ represents a permutation, and $\MAP(a)=p$ as desired.
\end{proof}

The next step is to show that $\MAP$ defines a bijection $\BIJ$ from the set of equivalence classes of permutations to weighted Motzkin paths, where two permutations are equivalent if they map to the same \emph{unweighted} Motzkin paths.
\begin{definition}
	 For an equivalence class $\mathbf{E}_{p}=\{\pi\in \mathcal{S}| \MAP(\pi)=p\}$ of permutations let $\BIJ(\mathbf{E}_{p})=p$, and let the weight of step $k$ be the sum of weights of node $k$ over permutations in $\mathbf{E_{p}}$.
\end{definition}

\begin{theorem}
The mapping $\BIJ(\mathbf{E})$ is a bijection from the set of equivalence classes of permutations (with the above definition of equivalent) to the set of weighted Motzkin paths, with weights
\begin{equation}\label{eq:weights 2} 
	N_{h}=E_{h}=\qxh{h+1} \mbox{ and }  S_{h}=F_{h}=\qh{h},
\end{equation}
such that the sum of weights of permutations in $\mathbf{E}$ is the weight of $\BIJ(\mathbf{E})$.
\end{theorem}
\begin{proof}
	Assume $\MAP$ maps node $k$ to an $E$ step at height $h$. Then the pair of left and right shapes is either $(\emptyset, \emptyset)$ or $(\rightarrow, \rightarrow)$. Further, to the left of node $k$ there must be $h+m$ shape pairs in the set $\{(\emptyset,\rightarrow),(\emptyset,\rightleftharpoons)\}$ (corresponding to $N$ steps) and $m$ shape pars in the set  $\{(\leftarrow,\emptyset),(\rightleftharpoons,\emptyset)\}$ (corresponding to $S$ steps). The nodes corresponding to $E$ and $F$ steps to the left of node $k$ may be disregarded in this discussion.
	
	Now, if the shape pair of node $k$ is $(\emptyset, \emptyset)$, there are $h$ arcs going over node $k$, and since all these arcs starts at a node to the left of node $k$, they are drawn after node $k$ is visited. Therefore node $k$ gets the weight $xq^{h}$. 
	
	If the shape pair is $(\rightarrow, \rightarrow)$ there is $h$ possibilities for the incoming arc (call this arc $A$). These give weights $1,\ldots, q^{h-1}$ depending on the number of arcs  with start node between the start node of arc $A$ and node $k$.
	
	Thus, in the image of the equivalence class, a step $E$ at height $h$ is given a total weight of $\qxh{h+1}$ as required.
	
	The cases of $F$, $N$ and $S$ steps are similar, and the details omitted. 
	
	That $\BIJ$ is a bijection follows at once from the fact that $\MAP$ is into the set of Motzkin paths, and $\BIJ$ is defined from the set of equivalence classes that maps to the same Motzkin paths.
\end{proof}

The step weights produced by $\BIJ$ are of the right form, but not exactly what we want. Let $\mathcal{M}_{n}^{\ast}$ denote the set of weighted Motzkin paths with weights given by \eqref{eq:weights 2}. A bijection $\BIJB$ from $\mathcal{M}_{n+1}^{\ast}$ to $\mathcal{M}_{n}$ will finally give paths with the correct weights.  
\begin{definition}
For $p$ in  $\mathcal{M}_{n+1}^{\ast}$ and for $k\in [n]$, if steps $k$ and $k+1$ is $x$ and $y$, let step $k$ in $\BIJB(p)$ be given by\\
\begin{center}
\begin{tabular}{c||c|c|c|c}$x\backslash y$ & $E$ & $F$ & $N$ & $S$ \\\hline\hline $E$ & $E$ & $S$ & $E$ & $S$ \\\hline $F$ & $N$ & $F$ & $N$ & $F$ \\\hline $N$ & $N$ & $F$ & $N$ & $F$ \\\hline $S$ & $E$ & $S$ & $E$ & $S$\end{tabular}
\end{center}
and have the same weight as step $k+1$ in $p$.
\end{definition}
\begin{theorem}
The mapping $\BIJB$ is a bijection from $\mathcal{M}_{n+1}^{\ast}$ to $\mathcal{M}_{n}$.
\end{theorem}
\begin{proof}[Proof (sketch)]
That $\BIJB$ give the correct step weights follows effortlessly from the definition. To show that $\BIJB$ is a bijection, the inverse mapping is easily derived. See \cite{Parv2006} for details.  
\end{proof}

%%%
%%%
%%%

\section{Closed forms}
Let $C(t)$ be the Catalan function, $C(t)=\frac{1-\sqrt{1-4t}}{2t}$. It is well known that $C(\gamma t)^{2}$ is the generating function for (bi-coloured) Motzkin paths in which each step have weight $\gamma$.

Define $\bar{\alpha}_{i}^{j}=\{\alpha_{i},\ldots,\alpha_{j}\}$ and $\bar{\beta}_{i}^{j}=\{\beta_{i},\ldots,\beta_{j}\}$. 
Let $g_{k}(\bar{\alpha}_{1}^{k}, \bar{\beta}_{1}^{k},\gamma; t)$ be the generating function for Motzkin paths in which weights are given by 
\begin{align*}
&N_{h}S_{h+1}=\beta_{h} \mbox{ for }  h\leq k,\\
&E_{h}+F_{h}=\alpha_{h} \mbox{ for }  h\leq k,\\
&N_{h}=S_{h+1}=E_{h}=F_{h}=\gamma \mbox{ for } h>k.
\end{align*}

Decomposing on the first return to the $x$-axis (where $E$ and $F$ steps counts as returns), we find that 
\[g_{1}(\alpha_{1},\beta_{1},\gamma; t)=1+(\alpha_{1}t+\beta_{1}t^{2} C(\gamma t)^{2})g_{1}(\alpha_{1},\beta_{1},\gamma; t),\]
and in general
\[g_{k}(\bar{\alpha}_{1}^{k}, \bar{\beta}_{1}^{k},\gamma; t)=1+(\alpha_{1}t+\beta_{1}t^{2}g_{k-1}(\bar{\alpha}_{2}^{k}, \bar{\beta}_{2}^{k},\gamma; t))g_{k}(\bar{\alpha}_{1}^{k}, \bar{\beta}_{1}^{k},\gamma; t).\]
Now, to find the number of permutations with $k$ occurrences of 2--13, we can count weighted Motzkin paths with all weights truncated at $q^{k}$. This is formalised in the following theorem. 
\begin{theorem}
For $i\leq k$,
\[\Phi_{i,j}(n)=[q^{i}x^{j}t^{n}]g_{k}(\{\qh{1}+\qxh{1},\ldots,\qh{k}+\qxh{k}\},\{\qh{1}\qxh{2},\ldots, \qh{k}\qxh{k+1}\},\qh{k};t).\] 
\end{theorem}
Let $G_{k}(x, t)=\sum_{m,n}x^{m}t^{n}\Phi_{k,m}(n)$. By iteratively calculating $g_{k}$ and differentiating with respect to $q$, we find that
\[
	G_{0}(x, t)=\frac{C(t)}{1 - t x C(t)},
\]
\[
	G_{1}(x, t)=\frac{C(t)(-1 + C(t)+ x(2 -  C(t)))(1 - C(t))^2}
	{(2 - C(t))(1 - x t C(t))^2}
\]
and
\begin{align*}
	G_{2}(x, t)=&\frac{2(1 - C(t))^3}{(2 - C(t))^{3}(1 - x t C(t))^3} 
	\Big(\\
	&-x^3(2 - C(t))^3(1 - C(t))\\
	&+ x^2(2 - C(t))^2 (3 - 8 C(t) + 4 C(t)^2)\\ 
	&- x (3 - 20 C(t) + 37 C(t)^2 - 24 C(t)^3 + 5 C(t)^4)	\\
	&-(1 - C(t)) (1 - 5  C(t) + 2  C(t)) 
        \Big).
\end{align*}

\subsection{Extracting coefficients}
The generating functions can be written in the form $P(C,x)(2-C)^{-a}(1- x t C)^{-b}$ for integers $a$ and $b$, and where $P(C,x)$ is a polynomial in $C$ and $x$. This allows for a routine, but lengthy, method for extracting coefficients. Consider as an example 
\[G_{1}(x,t)=\frac{A+xB}{(1 - x t C(t))^2}\]
where $A=\frac{C(t)(1 - C(t))^3}{C(t)-2}$ and $B=C(t)(1 - C(t))^2$. Expanding $G_{1}(x, t)$
in powers of $x$, we find that
\[[x^{k}] G_{1}(x, t)=(k+1)A t^{k}C(t)^{k}+k B t^{k-1}C(t)^{k-1}.\]
Noting that the above may be written as a sum of powers of $\sqrt{1-4t}$, coefficients may be extracted by applying the binomial theorem. See \cite{Parv2006} for details. 
\begin{theorem}
\begin{align*}
	\Phi_{0,m}(n)=&\binom{2n - m}{n - m}\frac{m + 1}{n+1},\\
	\Phi_{1,m}(n)=&\binom{2n - m}{n  - m-1}\frac{(m + 1) n^2 + 3(m - 1)n + m^2 + 9m + 2}{(n+2)(n+3)},\\
	\Phi_{2,m}(n)=&\binom{2n - m}{n - m - 1}\frac{1}{2(n + 2)(n + 3)(n + 4)}\\
	&\Big(
		(m + 1)n^4 - (6 - 5m + m^2)n^3 - (29 - 32m + 3m^2)n^2 \\&- (66 - 72m + 12m^2 + 2m^3)n-20 + 54m - 28m^2 - 6m^3
		\Big).
\end{align*}
\end{theorem}

%%%
%%%
%%%

\section{Other patterns}
There are 12 patterns of type (1,2) or (2,1). As shown by Claesson \cite{Clas2001}, these fall into three equivalence classes with respect to distribution of non-cyclic occurrences in permutations, namely
\begin{eqnarray*}
\mbox{\{1--23, 12--3, 3--21, 32--1\}, \{1--32, 21--3, 23--1, 3--12\} and \{13--2, 2--13, 2--31, 31--2\}.}
\end{eqnarray*}

It's only natural to ask about equivalence classes with respect to cyclic occurrences of patterns of type (1,2) and (2,1). Unfortunately, there are a lot of them. We conjecture that the 144 possible distributions fall into 106 equivalence classes. In any case 106 is lower bound. The conjectured classes of size 2 or more are given in Table 1.%\ref{table:equivalences}. 
\begin{conjecture}
The distributional relations in Table 1 %\ref{table:equivalences} 
holds, and the table includes all such relations.
\end{conjecture}
\begin{table}[htdp]\label{table:equivalences}
\begin{center}
\begin{tabular}{l}
(31--2, 31--2) \eqc (31--2, 2--31)\\% Reorder decent blocks\\ 
(13--2, 31--2) \eqc (13--2, 2--13)\\
(13--2, 13--2) \eqd (2--13, 2--13) \eqc (2--13, 31--2)\\
(2--31, 31--2) \eqc (2--31, 2--31)\\% Reorder decent blocks\\

(31--2, 3--21) \eqc (31--2, 32--1)\\
(2--31, 3--21) \eqc (2--31, 32--1)\\

(31--2, 3--12) \eqc (31--2, 23--1)\\
(13--2, 3--12) \eqc (13--2, 21--3)\\
(2--13, 3--12) \eqc (2--13, 21--3)\\
(2--31, 3--12) \eqc (2--31, 23--1)\\

(1--23, 31--2) \eqc (1--23, 2--13) \eqc(3--21, 31--2) \eqc (3--21, 2--31)\\
(3--21, 2--13) \eqc (1--23, 2--31)\\
(12--3, 31--2) \eqc (12--3, 2--13) \eqc (12--3, 2--31)\\
(32--1, 31--2) \eqc (32--1, 2--13) \eqc (32--1, 2--31)\\

(3--21, 13--2) \eqc (1--32, 13--2)\\

(1--32, 31--2) \eqc (1--32, 2--13) \eqc (1--32, 2--31)\\
(3--12, 31--2) \eqc (3--12, 2--13) \eqc (3--12, 2--31)\\
(21--3, 31--2) \eqc (21--3, 2--31)\\
(23--1, 31--2) \eqc (23--1, 2--13)\\

(1--23, 1--23) \eqd (12--3, 12--3)\\
(3--21, 3--21) \eqc (1--32, 32--1)\\
(1--32, 3--21) \eqc (3--21, 32--1)\\

(1--23, 3--12) \eqc (1--32, 21--3)\\
(3--21, 3--12) \eqc (1--32, 23--1)\\
(3--21, 21--3) \eqc (1--23, 23--1)\\
(1--32, 3--12) \eqc (1--23, 21--3) \eqc (3--21, 23--1)\\
(1--32, 1--32) \eqd (3--12, 3--12) \eqc (23--1, 23--1)
\end{tabular} 
\caption{Equivalences among cyclic occurrences of patterns of type (1,2) and (2,1).}
\end{center}
\end{table}

Let $\Pi(p_{b}, p_{w})$ denote the distribution of occurrences of $(p_{b}, p_{w})$ in permutations, where $(p_{b}, p_{c})$ means that we count occurrences of $p_{b}$ between cycles and of $p_{w}$ within cycles. Let $\Pi(p_{b}, p_{w}; C)$ denote the bivariate  distribution of cycles and occurrences of $(p_{b}, p_{w})$. 

Write $(p_{b}, p_{w})$ \eqd $(q_{b}, q_{w})$ if $\Pi(p_{b}, p_{w})=\Pi(q_{b}, q_{w})$, and $(p_{b}, p_{w})$ \eqc $(q_{b}, q_{w})$ if $\Pi(p_{b}, p_{w}; C)=\Pi(q_{b}, q_{w}; C)$.

First note that there are 8 \emph{diagonal} classes.
\begin{theorem}
The following distributional equivalences holds.
\begin{align*}
&\mbox{(13--2, 13--2) \eqd (2--13, 2--13),}\\ 
&\mbox{(1--23, 1--23) \eqd (12--3, 12--3) and}\\ 
&\mbox{(1--32, 1--32) \eqd (3--12, 3--12) \eqc (23--1, 23--1).}
\end{align*}
\end{theorem}
\begin{proof}
The ``$\sim$'' cases follow from Theorem \ref{t:onedim} and the non-cyclic equivalence classes.

It remains to show that (3--12, 3--12) \eqc (23--1, 23--1). Given a permutation $\pi$ in cycle form 
\[
C(\pi)=(c_{1}^{1} c_{2}^{1}\cdots c_{i_{1}}^{1})(c_{1}^{2} c_{2}^{2}\cdots c_{i_{2}}^{2})\cdots (c_{1}^{k} c_{2}^{k}\cdots c_{i_{k}}^{k}),
\]
let
\[
\hat{C}(\pi)=(d_{i_{k}}^{k} \cdots d_{1}^{k})\cdots (d_{i_{1}}^{1} \cdots d_{1}^{1}),
\]
where $d_{i}^{j}=n+1-c_{i}^{j}$. Write $\hat{C}(\pi)$ in standard cycle form. The result is a permutation, say $D(\pi)$, such that each occurrence, between or within cycles, of 3--12 in $C(\pi)$ corresponds exactly to an occurrence of 23--1 in $D(\pi)$. Furthermore, the cycle structure is obviously preserved.  \end{proof}

The seven patterns involved the above theorem share the property that they are equidistributed with the non-cyclic occurrences. Let $\Pi(p)$ denote the distribution of non-cyclic occurrences of the pattern $p$. 
\begin{theorem}\label{t:onedim}
We have
\begin{align*}
&\mbox{$\Pi$(2--13, 2--13) = $\Pi$(2--13),}\\
&\mbox{$\Pi$(13--2, 13--2) = $\Pi$(13--2),}\\
&\mbox{$\Pi$(1--23, 1--23) = $\Pi$(1--23),}\\
&\mbox{$\Pi$(12--3, 12--3) = $\Pi$(12--3),}\\
&\mbox{$\Pi$(1--32, 1--32) = $\Pi$(1--32),}\\
&\mbox{$\Pi$(23--1, 23--1) = $\Pi$(23--1) \mbox{ and }}\\
&\mbox{$\Pi$(3--12, 3--12) = $\Pi$(3--12).} 
\end{align*}
\end{theorem}
\begin{proof}
We use a standard bijection between permutations written in standard cycle form and permutations. Given a permutation $\pi$ in cycle form,
\[
C(\pi)=(c_{1}^{1} c_{2}^{1}\cdots c_{i_{1}}^{1})(c_{1}^{2} c_{2}^{2}\cdots c_{i_{2}}^{2})\cdots (c_{1}^{k} c_{2}^{k}\cdots c_{i_{k}}^{k}),
\]
map it to the permutation 
\[
\tilde \pi=c_{1}^{1} c_{2}^{1}\cdots c_{i_{1}}^{1}c_{1}^{2} c_{2}^{2}\cdots c_{i_{2}}^{2}\cdots c_{1}^{k} c_{2}^{k}\cdots c_{i_{k}}^{k}.
\]
Note that the bijection preserves the occurrences of each of the 7 patterns. This is true as we have the restrictions 
\[
c_{m}^{j}> c_{1}^{j+1}, m=1,\ldots, i_{j}, j=1, \ldots k-1.
\]
\end{proof}

%%%
%%%
%%%

\subsection{Increasing cycle order}
Using the standard cycle form and listing cycles in \emph{decreasing} order with respect to the cycles minimal elements is equivalent to listing cycles in \emph{increasing} order. 
\begin{theorem}  Let $\Pi^{d}(p_{b}, p_{w}; C)$ and $\Pi^{i}(p_{b}, p_{w}; C)$ denote the distribution of cyclic occurrence of some pattern pair $(p_{b}, p_{w})$ when the cycles are listed in decreasing respectively increasing order. 
Let $abc$ be a permutation of $[3]$. Then
\begin{align*}
&\Pi^{d}(\mbox{a--bc}, p_{w}; C)=\Pi^{i}(\mbox{bc--a}, p_{w}; C), \mbox{ and}\\
&\Pi^{d}(\mbox{ab--c}, p_{w}; C)=\Pi^{i}(\mbox{c-ab}, p_{w}; C).
\end{align*}
\end{theorem}
\begin{proof}
If $a$ is in a cycle to the left of a cycle containing $bc$ when the cycles are listed in decreasing order, it is to the right when the cycles are listed in increasing order.
\end{proof}
Writing cycles with the maximal element first also gives trivial equivalences. 
\begin{theorem}
Let $\hat\Pi^{x}$ denote the distributions when cycles are started with their maximal elements, and cycles are ordered in increasing $(x=i)$ or decreasing $(x=d)$ order. For a pattern $p$ let $r(p)$ denote the reverse pattern. Then
\begin{align*}
&\hat\Pi^{d}(p_{b}, p_{w})=\Pi^{i}(r(p_{b}), r(p_{w})), \mbox{ and}\\
&\hat\Pi^{i}(p_{b}, p_{w})=\Pi^{d}(r(p_{b}), r(p_{w})).
\end{align*} 

\end{theorem}

%%%
%%%
%%%

%\section{What about $y$, and other questions}
%Disregarding $q$ we get a continued fraction form of the ordinary generating function for the unsigned Stirling numbers of the first kind, $|S(n,k)|$. Although probably known, the author have not been able to find this in the literature.
%\begin{corollary}
%\begin{align*}
%	1+x t F(1, x, 1, t)&=1+\cfrac{x t}{1-(1+x)t
%   	-\cfrac{1(1+x)t^{2}}{1- (3+x)t
%	-\cfrac{2(2+x)t^{2}}{1- (5+x)t \dotsb}}}\\
%	&=
%	\cfrac{1}{1-
%		\cfrac{x t}{1-
%			\cfrac{t}{1-
%				\cfrac{(1+x)t}{1-
%					\cfrac{2t}{1-
%				\cfrac{(2+x) t}{1-
%					\cfrac{3t}{1-\cdots}}}}}}}\\
%	&=
%	\sum_{n=0}^{\infty}\sum_{k=0}^{\infty}|S(n,k)|x^{k}z^{n}.
%\end{align*}
%\end{corollary}

%%%%
%%%%
%%%%

\section{What about $y$?}
Expanding $F(1,x,y,t)$ we are quickly led to conjecture that $F(1,x,y,t)$ is the generating function for a product of Stirling numbers and binomial coefficients. Using the same bijection as in the proof of Theorem \ref{theorem:main}, we can prove this. 
\begin{theorem}
	\[[x^{i}y^{j}t^{n}]F(1,x,y,t)=\binom{i+j}{j}|S(n, i+j)|.\]
%	\[F(1,x,y,t)=\sum_{i}\sum_{j}\sum_{n}\binom{i+j}{j}|S(n, i+j)|x^{i}y^{j}t^{n}.\]
\end{theorem}
In other words, $[x^{i}y^{j}z^{n}]F(1,x,y,t)$ is the number of permutations of $[n]$ with $i+j$ cycles of which $i$ are marked. We will call these \emph{marked permutations}, and denote the set of marked permutations of length $n$ with $\underline{\mathcal{S}}_{n}$. In Table 2 the elements of $\underline{\mathcal{S}}_{3}$ are listed. As the proof is much the same as that of Theorem \ref{theorem:main}, we only sketch it here.

\begin{table}\label{table:marked}
\begin{center}
\renewcommand\arraystretch{1.5}
\begin{tabular}{c||c|c|c|c}
&0&1&2&3\\
\hline\hline
&&(123)  & (23)(1) & (3)(2)(1)\\
0&&(132) & (3)(12) &\\
& & &  (2)(13) &\\
\hline
&$\underline{(123)}$ & $\underline{(23)}(1), (23)\underline{(1)}$ &$\underline{(1)}(2)(3)$\\
1&$\underline{(132)}$ & $\underline{(2)}(13) , (2)\underline{(13)}$ &$(1)\underline{(2)}(3)$\\
     &                             & $\underline{(3)}(12), (3)\underline{(12)}$ & $(1)(2)\underline{(3)}$\\
\hline
&$\underline{(1)}\underline{(23)}$ & $\underline{(3)}\underline{(2)}(1)$&\\
2&$\underline{(2)(13)}$                   & $(3)\underline{(2)(1)}$&\\
&$\underline{(3)(12)}$                   & $\underline{(3)}(2)\underline{(1)}$&\\
\hline
3&$\underline{(3)(2)(1)}$&&
\end{tabular}
\caption{The set of marked permutations of length 3. The marked cycles are \underline{underlined}.}
\end{center}
\end{table}
\begin{proof}[Proof (sketch)]
We again use the arc representation. Give node $k$ weight $x$ if it is the first element in an unmarked cycle, and weight $y$ if it is the first in a marked cycle. 

Reasoning as in the proof of Theorem  \ref{theorem:main} shows that $\MAP$
defines a bijection from equivalence classes of permutations with the above weighting to weighted Motzkin paths with weights
\[N_{h}=E_{h}=h+x+y, S_{h}=F_{h}=h.\]
The result follows after application of \cite[Theorem 1]{Flaj1980}.
\end{proof}

%%%
%%%
%%%

\subsection{What about $q$ and $y$?}
In light of the above, $F(q,x,y,t)$ should count the number of permutations with respect to length, cycles, marked cycles and occurrences of 2--13. Unfortunately, life is not that easy. For instance, $[t^{3}]F(q,1,1,t)=14+8q+q^{2}$, but in the set of 24 marked permutations of length 3 there are only two single occurrences of 2--13. 

Perhaps marked permutations are not the natural object for studying $F(q,x,y,t)$. As the number of marked permutations of length $n$ is $(n+1)!$, we should look for a nice (weight preserving) bijection between  $\underline{\mathcal{S}}_{n}$ and $\mathcal{S}_{n+1}$. So far, we have not found such a bijection.

\bibliographystyle{abbrv}
\bibliography{/Users/robertp/Documents/Projects/biblio.bib}
\end{document}